\documentclass[11pt]{article}
\usepackage{latexsym}
\usepackage{graphicx}
\usepackage{amsmath}
\usepackage{amssymb}
\usepackage{color}
\usepackage{amscd}
\usepackage{epsfig}
\usepackage{cite}
\newtheorem{problem}{Problem}[section]

\newtheorem{lemma}[problem]{Lemma}
\newtheorem{theorem}[problem]{Theorem}
\newtheorem{prop}[problem]{Proposition}
\newtheorem{corollary}[problem]{Corollary}

\title{The $L$-functions of twisted Witt extensions}
\author{Chunlei Liu
\\
\small {School of Mathematical Science, Beijing Normal University,
Beijing 100875. {\it Email:} clliu@bnu.edu.cn}}
\begin{document}
\maketitle
{\bf Abstract. }The $L$-function of a non-degenerate twisted Witt
extension is proved to be a polynomial. Its Newton polygon is
proved to lie above the Hodge polygon of that extension. And the
Newton polygons of the Gauss-Heilbronn sums are explicitly
determined, generalizing the Stickelberger theorem.


\section{Introduction}
\hskip 0.2in Let $p$ be a prime number,
$\mathbb{Z}_p:=\lim\limits_{\leftarrow}\mathbb{Z}/(p^n)$ be the
ring of $p$-adic integers, and $\mathbb{Q}_p$ its fraction field.
For every positive integer $n$, denote by $\mu_n$ the group of
$n$-th roots of unity. Fix an algebraic closure
$\overline{\mathbb{Q}}_p$ of $\mathbb{Q}_p$. Then, for every
positive power $q$ of $p$,
$\mathbb{F}_p:=\mathbb{Z}_p[\mu_{q-1}]/(p)$ is a finite field of
$q$ elements. Every $a$ in $\mathbb{F}_q$ has a unique
representative $\omega(a)$ in $\mu_{q-1}\cup\{0\}$. It is known as
the Teichm\"{u}ller representative.

Let $m$ be a positive integer. For every ring $A$, the ring
$W_m(A)$ of Witt vectors of length $m$ with coefficients in $A$ is
the set $A^m$ endowed with the Witt ring structure. It is known
that $W_m(\mathbb{F}_q)$ is isomorphic to
$\mathbb{Z}_p[\mu_{q-1}]/(p^m)$. The isomorphism is given by
$$\iota:(a_0,\cdots,a_{m-1})\mapsto\sum\limits_{i=0}^{m-1}\omega(a_i^{p^{-i}})p^i.$$

Fix a character $\Psi_0$ of $\mathbb{Z}_p/(p^m)$ of exact order
$p^m$. Then
$\Psi_q:=\Psi_0\circ\iota\circ\text{Tr}_{W_m(\mathbb{F}_q)/W_m(\mathbb{F}_p)}$
is a character of $W_m(\mathbb{F}_q)$. Let $f$ be a Witt vector of
length $m$ with coefficients in
$\mathbb{F}_q[x_1^{\pm1},\cdots,x_n^{\pm1}]$ and $\chi$ a
character of $(\mathbb{F}_q^{\times})^n$. We introduce the
exponential sum
$$S(f,\chi):=(-1)^{n-1}\sum\limits_{x\in(\mathbb{F}_q^{\times})^n}
\chi(x)\Psi_q(f(x)).$$ If $\chi$ is trivial, it is the exponential
sum studied by Liu-Wei [LW]. If $m=1$, it is the twisted
exponential sums studied by Adolphson-Sperber [AS, AS2]. In this
paper, we assume that the first coordinate of $f$ is nonzero. That
obviously loses no generality.

Consider the extension
$$\mathbb{F}_q(x_1,\cdots,x_n)\hookrightarrow
\mathbb{F}_q(x_1^{\frac{1}{q-1}},\cdots,x_n^{\frac{1}{q-1}})(y_0,\cdots,y_{m-1}),$$
where $(y_0,\cdots,y_{m-1})$ satisfies the equation of Witt
vectors $$(y_0^p,\cdots,y_{m-1}^p)-(y_0,\cdots,y_{m-1})=f(x).$$ It
is called the twisted Witt extension associated to $f$. Its Galois
group $G$ is isomorphic to
$W_m(\mathbb{F}_p)\times(\mathbb{F}_q^{\times})^n$. The
isomorphism is given by
$$\kappa:g\mapsto (g(y_0,\cdots,y_{m-1})-(y_0,\cdots,y_{m-1}),\frac{g(x^{\frac{1}{q-1}})}{x^{\frac{1}{q-1}}}),$$
where $$\frac{g(x^{\frac{1}{q-1}})}{x^{\frac{1}{q-1}}}
=(\frac{g(x_1^{\frac{1}{q-1}})}{x_1^{\frac{1}{q-1}}},\cdots,
\frac{g(x_n^{\frac{1}{q-1}})}{x_n^{\frac{1}{q-1}}}).$$Write
$\mathbb{F}:=\lim\limits_{\rightarrow}\mathbb{F}_{p^k}$. It is an
algebraic closure of $\mathbb{F}_p$. Let $\tilde{x}$ be a class of
degree $k$ in
$(\mathbb{F}^{\times})^n/\text{Gal}(\mathbb{F}/\mathbb{F}_q)$, and
$\text{Fr}_{\tilde{x}}$ the Frobenius element of $G$ at
$\tilde{x}$. Then we can show that
$$\kappa(\text{Fr}_{\tilde{x}})=(\text{Tr}_{W_m(\mathbb{F}_{q^k})/W_m(\mathbb{F}_p)}(f(\bar{x})),
\text{N}_{\mathbb{F}_{q^k}/\mathbb{F}_q}(\bar{x})),$$ where
$\bar{x}$ is an element of $(\mathbb{F}_{q^k}^{\times})^n$
representing $\tilde{x}$. Note that
$\rho:=(\Psi_p,\chi)\circ\kappa$ is a character of $G$. The Artin
$L$-function of $\mathbb{F}_q(x_1,\cdots,x_n)$ associated to
$\rho$ is
$$L_{f,\chi}(t)=\prod\limits_{x\in(\mathbb{F}^{\times})^n/\text{Gal}(\mathbb{F}/\mathbb{F}_q)}
(1-\rho(\text{Fr}_{\tilde{x}})t^{\deg(x)})^{(-1)^n}.$$ It is an
analytic expression for the arithmetic of the twisted Witt
extension associated to $f$. One can show that
$$L_{f,\chi}(t)=
\exp(\sum\limits_{k=1}^{+\infty} S_k\frac{t^k}{k}),$$where
$S_k:=S(f,\chi\circ\text{N}_{\mathbb{F}_{q^k}/\mathbb{F}_q})$.

For every ring $A$, define $$[\cdot]:A\rightarrow W_m(A),\
a\mapsto(a,0,\cdots,0),$$ and
$$V:A\rightarrow W_m(A),\ (a_0,\cdots,a_{m-1})\mapsto(0,a_0,\cdots,a_{m-2}).$$
According to [LW, \S 0], the Witt vector $f$ has a unique
decomposition of the form$$f=\sum
\limits_{i=0}^{m-1}\sum\limits_{u\in I_i}V^i([a_{iu}x^u]),\
I_i\subseteq\mathbb{Z}^n,\ a_{iu}\in\mathbb{F}_q^{\times}.$$ For
each $i=0,\cdots,m-1$, we denote by $\Delta_i$ the convex hull in
$\mathbb{Q}^n$ of $I_i$ and the origin. We define $\Delta$ to be
the $m$-tuple $(\Delta_0,\cdots,\Delta_{m-1})$, and denote by
$\Delta_{\infty}$ the convex hull in $\mathbb{Q}^n$ of
$\cup_{i=0}^{m-1}p^{m-i-1}\Delta_i$.

In this paper, we assume that $\Delta_{\infty}$ generates
$\mathbb{Q}^n$. Suppose that $\Delta_{\infty}$ generates a
subspace of dimension $l$. In that $l$-dimensional subspace,
choose $l$ linearly independent integral vectors
$(\alpha_{i1},\cdots,\alpha_{in})$ ($i=1,\cdots,l$) that span a
parallelotope of the smallest volume. Then there are integral
matrices $P$ and $Q$ with determinants $\pm1$ such that
$P(\alpha_{ij})_{1\leq i\leq l,1\leq j\leq n}=(I,0)Q$. It follows
that we can choose integral vectors
$(\alpha_{i1},\cdots,\alpha_{in})$ ($i=l,\cdots,n$) such that
$(\alpha_{ij})_{1\leq i\leq n,1\leq j\leq n}$ has determinant
$\pm1$. Making the change of variables
$$\left\{\begin{array}{l}
  y_1=x_1^{\alpha_{11}}\cdots x_n^{\alpha_{1n}} \\
  \vdots \\
  y_l=x_n^{\alpha_{l1}}\cdots x_n^{\alpha_{ln}}, \\
  z_1=x_1^{\alpha_{l+1,1}}\cdots x_n^{\alpha_{l+1,n}} \\
  \vdots \\
  z_{n-l}=x_n^{\alpha_{n1}}\cdots x_n^{\alpha_{nn}}, \\
\end{array}\right.$$ we see that
$$S(f,\chi)=(-1)^{n-1}\sum\limits_{z\in(\mathbb{F}_q^{\times})^{n-l}}
\chi_2(z)\sum\limits_{y\in(\mathbb{F}_q^{\times})^l}
\chi_1(y)\Psi_q(g(y))$$ for some character $(\chi_1,\chi_2)$ of $
(\mathbb{F}_q^{\times})^n$, and some Witt vector $g$ of the form
$$g=\sum \limits_{i=0}^{m-1}\sum\limits_{v\in
J_i}V^i([b_{iv}y^v]),\ J_i\subseteq\mathbb{Z}^l,\
b_{iv}\in\mathbb{F}_q^{\times}.$$ Therefore assuming that
$\Delta_{\infty}$ generates $\mathbb{Q}^n$ loses no generality.

We call $f$ non-degenerate with respect to $\Delta$ if for every
face $\tau$ of $\Delta_{\infty}$ that does not contain 0, the
system $\overline{_1f}^{\tau} =\cdots=\overline{_nf}^{\tau}=0$ has
no common solution in $(\mathbb{F}^{\times})^n$, where
$$\overline{_jf}^{\tau}=\sum
\limits_{i=0}^{m-1}\sum\limits_{p^{m-i-1}u\in\tau}u_ja_{iu}^{p^{m-i-1}}x^{p^{m-i-1}u}.$$In
this paper we assume that $f$ is non-degenerate. Generalizing
corresponding results of Adolphson-Sperber [AS2, Corollary 2.12]
and Liu-Wei [LW, Theorem 1.3], we establish the following theorem.
\begin{theorem}For each $\chi$, $L_{f,\chi}$ is a polynomial.
\end{theorem}

Since $L_{f,\chi}$ is a polynomial with coefficients in
$\mathbb{Q}(\mu_{(q-1)p^m})$, it is interesting to know the prime
decomposition of its reciprocal roots. Fix an embedding of
$\overline{\mathbb{Q}}$ into $\overline{\mathbb{Q}}_p$. Knowing
the prime decomposition of all reciprocal roots is equivalent to
knowing the Newton polygon of $L_{f,\chi}$. The Newton polygon of
a polynomial $1+\sum\limits_{i=0}^k\alpha_it^i$, with respect to a
valuation $\text{ord}$ of $\overline{\mathbb{Q}}_p$, is the convex
hull in $\mathbb{Q}^2$ of the points $(0,0)$ and
$(i,\text{ord}(\alpha_i))$ ($i=1,\cdots,k$). Denote by
$C(\Delta_{\infty})$ the cone in $\mathbb{Q}^n$ generated by
$\Delta$. There is a $\mathbb{R}_{\geq0}$-linear degree function
$u\mapsto\deg(u)$ on $C(\Delta_{\infty})$ such that $\deg(u)=1$
when $u$ lies on a face of $\Delta_{\infty}$ that does not contain
the origin. For every integer $s$, write
$L_s(\Delta_{\infty}):=C(\Delta_{\infty})\cap(\frac{s}{q-1}+\mathbb{Z}^n)$.
There is a least positive integer $D$ such that $D\deg
L_0(\Delta_{\infty})\subset\mathbb{Z}$. For every natural number
$k$, we denote by $W_s(k)$ the number of points of degree
$\frac{k}{D(q-1)}$ in $L_s(\Delta_{\infty})$. Define
$P_s(t)=(1-t^{D(q-1)})^n\sum\limits_{k=0}^{+\infty}W_s(k)t^k$. By
Corollary 3.9, $P_s(t)$ is a polynomial with coefficients in the
set of natural numbers. The degree-$M$ Hodge polygon of a
polynomial $\sum\limits_{i=0}^l\alpha_it^i$ with non-negative
coefficients is the polygon in $\mathbb{Q}^2$ with vertices at the
points $(0,0)$ and
$(\sum\limits_{i=0}^k\alpha_i,\sum\limits_{i=0}^k\frac{i}{M}\alpha_i)$
($i=0,\cdots,l$). From now on, we write $q=p^a$ and denote by $s$
the integer such that, for each $x\in\mathbb{F}_q^{\times}$,
$\chi(x)=\omega(x)^{-s}$ in $\overline{\mathbb{Q}}_p$.
Generalizing corresponding results of Adolphson-Sperber [AS2,
Theorem 3.17] and Liu-Wei [LW, Theorem 1.3], we establish the
following theorem.
\begin{theorem}The Newton polygon of $L_{f,\chi}$ with respect to $\text{ord}_q$ lies
above the degree-$D(q-1)$ Hodge polygon of
$\frac{1}{a}\sum\limits_{i=0}^{a-1}P_{sp^i}(t)$ with the same
endpoints. In particular, it is of degree
$n!\text{Vol}(\Delta_{\infty})$.
\end{theorem}

We call $S(f,\chi)$ a Gauss-Heilbronn sum if $n=1$ and
$f=\sum\limits_{i=0}^{m-1}V^i([c_ix])$. Without loss of
generality, we assume that $c_0=1$. If $m=1$, the Gauss-Heilbronn
sum becomes a Gauss sum. And if $\chi$ is trivial, it is a
Heilbronn sum. Write
$\frac{s}{q-1}=-\sum\limits_{l=0}^{+\infty}s_lp^l$. Generalizing
the Stickelberger theorem, we establish the following theorem.
\begin{theorem}If $n=1$ and
$f=\sum\limits_{i=0}^{m-1}V^i([c_ix])$ with $c_0=1$, then the
Newton polygon of $L_{\chi}(t)$ with respect to $\text{ord}_q$
coincides with the polygon with vertices at $(0,0)$ and the points
$$(k,\frac{(k-1)k}{2p^{m-1}}+\frac{k(s_0+\cdots+s_{a-1})}{ap^{m-1}(p-1)}),\ k=1,\cdots,p^{m-1}.$$
\end{theorem}

Theorem 1.3 is the main contribution of this paper. In the proof
of Theorem 1.3, we choose a suitable basis of the space in
question so that the corresponding matrix of the Dwork operator is
amenable, and then we relate that matrix to a matrix from the
coefficients of the exponential function.

{\bf Acknowledgement. }This work is supported by NSFC Grant No.
10371132, by Project 985 of Beijing Normal University, by the
Foundation of Henan Province for Outstanding Youth, and by the
Morningside Center of Mathematics in Beijing. The author thanks
Daqing Wan for encouragement in the study of Gauss-Heilbronn sums.
\section{The $p$-adic trace formula}
We prove Theorem 1.1 after establishing a $p$-adic trace formula
relating $L_{f,\chi}(t)$ to the characteristic polynomials of an
operator on $p$-adic spaces.

Let
$$E(t)=\exp(\sum_{i=0}^{\infty}\frac{t^{p^i}}{p^i}) \in
\mathbb{Z}_p[[t]]$$ be the Artin-Hasse exponential series.
\begin{lemma}([LW, Lemma 2.2]) Let $l$ be a positive integer.
Then $\pi\mapsto E(\pi)$ is a bijection from the set of roots of
$\sum\limits_{i=0}^{\infty}\frac{t^{p^i}}{p^i}=0$ in
$\overline{\mathbb{Q}}_p$ with order $\frac{1}{p^{l-1}(p-1)}$ to
the set of all primitive $p^l$-th roots of unity in
$\overline{\mathbb{Q}}_p$.
\end{lemma}

By that lemma, we may choose, for each $l=1,\cdots,m$, a unique
root $\pi_l$ of $\sum\limits_{i=0}^{\infty}\frac{t^{p^i}}{p^i}=0$
in $\overline{\mathbb{Q}}_p$ with order $\frac{1}{p^{l-1}(p-1)}$
such that $E(\pi_l)=\Psi_0(1)^{p^{m-l}}$. Let $\pi$ be a
$D(q-1)$-th root of $\pi_m$, and $O:=\mathbb{Z}_p[\mu_{q-1},\pi]$.
For every $b\geq0$, write
$$L_s(b)=\{\sum_{u\in L_s(\Delta_{\infty})}a_ux^u :
 a_u\in  O, \text{ ord}_p(a_u)\geq b\deg(u)\}.$$
Define $$E_f(x) =\prod\limits_{i=0}^{m-1}\prod\limits_{u\in
I_i}E(\pi_{m-i}\omega(a_{iu})x^u).$$
\begin{lemma}([LW, Lemma 2.5]) We have $E_f\in L_0(\frac{1}{p-1})$.\end{lemma}

Let $\sigma$ be the element of
$\text{Gal}(\mathbb{Q}_p[\mu_{q-1},\pi]/\mathbb{Q}_p)$ that fixes
$\pi$ and becomes the $p$-power map when restricted to
$\mu_{q-1}$. Let $\sigma$ acts on $L(b)$ coefficientwise. For
every positive integer $k$, write
$$E_{f,q^k}(x)=\prod\limits_{l=0}^{ak-1}E_f^{\sigma^l}(x^{p^l}).$$
It is easy to see that
$$E_{f,q^k}(x)=\prod\limits_{j=0}^{k-1}E_{f,q}(x^{q^j}).$$
And, by the above lemma, we have $E_{f,q^k}\in
L_0(\frac{1}{p^{ak-1}(p-1)})$.
\begin{lemma}([LW, Lemma 2.6]) For
every positive integer $k$, we have
$$\Psi_{q^k}(f(x))
=E_{f,q^k}(\omega(x)),\ x\in\mathbb{F}_{q^k}.$$\end{lemma}
\begin{corollary}For every positive integer $k$, we have
$$S_k
=(-1)^{n-1}\sum\limits_{x\in\mu_{q^k-1}^n}x^{-d(1+q+\cdots+q^k)}E_{f,q^k}(x).$$
\end{corollary}

Write
$L_s(\Delta_{\infty})=C(\Delta_{\infty})\cap(\frac{d}{q-1}+\mathbb{Z}^n)$.
Then $L_0(\Delta_{\infty})+L_s(\Delta_{\infty})\subseteq
L_s(\Delta_{\infty})$.
 \begin{corollary}Let $k$ be a positive integer. If
$E_{f,q^k}(x)=\sum\limits_{u\in L_0(\Delta_{\infty})}a_ux^u $,
then
$$S_k
=(-1)^{n-1}(q^k-1)^n\sum\limits_{u\in L_s(\Delta_{\infty})}
a_{(q^k-1)u}.$$
\end{corollary}

 Define $$B_s:=\{\sum\limits_{u\in L_s(\Delta_{\infty})}a_ux^u\in
 L_s(\frac{1}{p-1}):
 \text{ ord}_p(a_u)-\frac{\deg(u)}{p-1}\rightarrow+\infty
 \text{ if }\deg(u)\rightarrow+\infty\}.$$
Then $B_s$ is a $B_0$-module.
 For every $\sum\limits_{u\in L_s(\Delta_{\infty})}a_ux^u\in B_s$, define
$$
\|\sum\limits_{u\in
L_s(\Delta_{\infty})}a_ux^u\|=\max\limits_{u\in
L_s(\Delta_{\infty})}p^{\frac{\deg(u)}{p-1}}|a_u|_p.$$ The
following lemma is easy.
\begin{lemma}The map $\|\cdot\|$ is a norm on $B_s$ over $O$, and
$B_s$ is complete with respect to that norm.\end{lemma}

For every $b\geq0$, write
$$L_s(b)=\{\sum\limits_{u\in L_s(\Delta_{\infty})}a_ux^u :
 a_u\in  O, \text{ ord}_p(a_u)\geq b\deg(u)\}.$$
 Then $L_s(b)$ is a $L_0(b)$-module.
 Define $$\psi_p:L_s(b)\rightarrow L_{sp^{a-1}}(pb),\
 \sum\limits_{u\in
L_s} a_ux^u\mapsto\sum\limits_{u\in L_{sp^{a-1}}} a_{pu}x^u.$$
\begin{lemma}The map $\phi_p:=\sigma^{-1}\circ\psi_p\circ E_f$
sends $B_s$ to $B_{sp^{a-1}}$. We call it a Dwork operator.
\end{lemma}
{\it Proof. }Let $g(x)\in B_s$. Then $gE_f\in L_s(\frac{1}{p-1})$.
So $\psi_p(gE_f)\in L_{sp^{a-1}}(\frac{p}{p-1})\subseteq
B_{sp^{a-1}}$. The lemma is proved.
\begin{corollary}The map $\phi_p^{ak}=\psi_p^{ak}\circ E_{f,q^k}$. It acts
on $B_s$, and is $O$-linear.
\end{corollary}

For every sequence $\{a_u\}\in O^{L_s(\Delta_{\infty})}$, we
define
$$\|\{a_u\}\|:=\max_{u\in L_s(\Delta_{\infty})}|a_u|_p.$$
\begin{lemma}The map $\|\cdot\|$ is a norm on $ O^{L_s(\Delta_{\infty})}$
over $ O$, and $ O^{L_s(\Delta_{\infty})}$ is complete with
respect to that norm.\end{lemma}

For each $v\in L_s(\Delta_{\infty})$, define a column vector
$(c_{uv})_{u\in L_{sp^{a-1}}(\Delta_{\infty})}$ with coefficients
in $O$ by
$$\phi_p(\pi^{p^{m-1}D(q-1)\deg(v)}x^v)
=\sum\limits_{u\in
L_{sp^{a-1}}(\Delta_{\infty})}c_{uv}\pi^{p^{m-1}D(q-1)\deg(u)}x^u.$$
\begin{lemma}We have
$\|(c_{uv})_{v\in L_s(\Delta_{\infty})}\|\rightarrow0$ if
$\deg(u)\rightarrow+\infty$.\end{lemma} {\it Proof. } Write
$$E_f(x)
=\sum\limits_{u\in L(\Delta_{\infty})}a_ux^u,\
\text{ord}_p(a_u)\geq\frac{\deg(u)}{p-1}.$$
Then$$\phi_p(\pi^{p^{m-1}D(q-1)\deg(v)}x^v)=\sum\limits_{u\in
L_{sp^{a-1}}(\Delta_{\infty})}a_{pu-v}\pi^{p^{m-1}D(q-1)(\deg(v)-\deg(u))}\pi^{p^{m-1}D(q-1)\deg(u)}x^u.$$
The lemma then follows from the fact that
$$\|\{a_{pu-v}\pi^{p^{m-1}D(q-1)(\deg(v)-\deg(u))}\}_{v\in L_s(\Delta_{\infty})}\|
\leq p^{-\deg(u)}.$$

\begin{corollary}The operator
$\phi_p$ is completely continuous.
\end{corollary}

By Serre [Se], for every positive integer $k$, the trace of
$\phi_p^{ak}$ over $O$ is well-defined and is equal to the trace
of the matrix of $\phi_p^{ak}$ with respect to any orthonormal
basis of $B_s$.
\begin{lemma}For
every positive integer $k$,
$$S_k
=(-1)^{n-1}(q^k-1)^n\text{Tr}_{B_s/ O}(\phi_p^{ak}).$$\end{lemma}
{\it Proof. }Let $g(x)\in B_s$. We have
$$\phi_p^{ak}(g)=\psi_p^{ak}(gE_{f,q^k}).$$Write
$E_{f,q^k}(x)=\sum\limits_{u\in L(\Delta_{\infty})}a_ux^u $. Then
$$\phi_p^{ak}(\pi^{p^{m-1}D(q-1)\deg(v)}x^v)=\sum\limits_{u\in
L_s(\Delta_{\infty})}a_{q^ku-v}\pi^{p^{m-1}D(q-1)\deg(v)}x^u.$$ So
the trace of $\phi_p^{ak}$ on $B_s$ over $O$ equals
$\sum\limits_{u\in L_s(\Delta_{\infty})} a_{(q^k-1)u}$. The lemma
then follows from Corollary 2.5.

Let $e_1=(1,0,\cdots,0)$, $\cdots$, $e_n=(0,\cdots,0,1)$. For each
$l=0,1,\cdots,n$, write
$$K_{s,l}=\bigoplus\limits_{1\leq
i_1<\cdots<i_l\leq n} B_se_{i_1}\wedge\cdots\wedge e_{i_l}$$ and
define
$$\phi_{p,l}:K_{s,l}\rightarrow K_{sp^{a-1},l},
\ ge_{i_1}\wedge\cdots\wedge e_{i_l}\mapsto
p^l\phi_q(g)e_{i_1}\wedge\cdots\wedge e_{i_l}.$$By Lemma 2.12, we
have the following chain-level trace formula.
\begin{prop}For every positive integer $k$,
$$S_k
=\sum\limits_{l=0}^n (-1)^{l+1} \text{Tr}_{K_{s,l}/
O}(\phi_{p,l}^{ak}).$$
\end{prop}

Define
$$\hat{E}_f(x):=\prod\limits_{j=0}^{+\infty}E_{f,q}(x^{q^j}).$$ And write
$$d\log\widehat{E}_f(x)=
\sum\limits_{k=1}^n\widehat{_kf}\frac{dx_k}{x_k}.$$
\begin{lemma}([LW, Corollary 3.8]) For
$k=1,\cdots,n$, we have $\widehat{_kf}\in B_0$, and
$$\widehat{_kf}\equiv
\sum\limits_{i=0}^{m-1}\sum\limits_{j=0}^{m-i-1}
\sum\limits_{\deg(p^{m-i-1}u)=1}u_k\omega(a_{iu}^{p^j})\pi^{D(q-1)p^{m-1}\deg(p^ju)}x^{p^ju}\
(\text{mod }\pi B_0).$$
\end{lemma}

By that lemma, $\hat{\partial}_j:g\mapsto (x_j\frac{\partial
}{\partial x_j}+\widehat{_jf})g$, $j=1,\cdots,n$, operate on
$B_s$. Obviously, they commute with each other. So, for each
$l=1,\cdots,n$,
$$\hat{\partial}:K_{s,l}\rightarrow K_{s,l-1},\
ge_{i_1}\wedge\cdots\wedge
e_{i_l}\mapsto\sum\limits_{k=1}^l(-1)^{k-1}\hat{\partial}_{i_k}(g)e_{i_1}\wedge\cdots\wedge
\hat{e}_{i_k}\wedge\cdots\wedge e_{i_l},\ i_1<\cdots<i_l$$ is
well-defined. We have $\hat{\partial}^2=0$. Thus we get a complex
$$K_{s,n}\stackrel{\hat{\partial}}{\rightarrow}K_{s,n-1}
\stackrel{\hat{\partial}}{\rightarrow}\cdots\stackrel{\hat{\partial}}{\rightarrow}K_{s,0}.$$
It is easy to check that
$\phi_{p,l-1}\circ\hat{\partial}=\hat{\partial}\circ \phi_{p,l}$.
That is, $\phi_p:=(\phi_{p,n},\cdots,\phi_{p,0})$ sends the
complex $(K_{s,\bullet},\hat{\partial})$ to the complex
$(K_{sp^{a-1},\bullet},\hat{\partial})$. So $\phi_p^{ak}$ operates
$O$-linearly on the complex $(K_{s,\bullet},\hat{\partial})$.
Therefore we have the following homological trace formula.
\begin{prop}For every positive integer $k$,
$$S_k=\sum\limits_{l=0}^n (-1)^{l+1}
\text{Tr}_{H_l(K_{s,\bullet},\hat{\partial})/ O}(\phi_p^{ak}).$$
Equivalently, $$L_{f,\chi}(t)=\prod\limits_{l=0}^n
\text{det}_O(1-\phi_p^at\mid
H_l(K_{s,\bullet},\hat{\partial}))^{(-1)^l}.$$
\end{prop}
Define $B=\oplus_{s=0}^{q-2}B_s$. For each $l=0,1,\cdots,n$, write
$$K_l=\oplus_{s=0}^{q-2}K_{s,l}=\bigoplus\limits_{1\leq
i_1<\cdots<i_l\leq n} Be_{i_1}\wedge\cdots\wedge e_{i_l}.$$
\begin{lemma}([LW, Proposition 6.1]) The $O$-module
$H_l(K_{\bullet},\hat{\partial})$ is 0 if $l=1,\cdots,n$, and is
free of finite rank if $l=0$.
\end{lemma}Form that lemma and the fact that
$$H_l(K_{\bullet},\hat{\partial})=\oplus_{s=0}^{q-2}H_l(K_{s,\bullet},\hat{\partial}),$$
we deduce the following corollary.
\begin{corollary}The $O$-module
$H_l(K_{s,\bullet},\hat{\partial})$ is 0 if $l=1,\cdots,n$, and is
free of finite rank if $l=0$.
\end{corollary}From that corollary and the homological trace formula, we deduce Theorem 1.1.
More precisely, we have the following
corollary.
\begin{corollary}The $O$-module
$H_0(K_{s,\bullet},\hat{\partial})$ is free of finite rank and
$$L_{f,\chi}(t)=\text{det}_O(1-\phi_p^at\mid
H_0(K_{s,\bullet},\hat{\partial})).$$
\end{corollary}
\section{Bases represented by homogenous elements}
We prove the following proposition.\begin{prop}The $O$-module
$H_0(K_{s,\bullet},\hat{\partial})$ is free of finite rank
$n!\text{Vol}(\Delta_{\infty})$. Moreover, it has a basis
represented a set $V_s$ of homogenous elements such that
$$P_s(t)=\sum\limits_{k=0}t^k\sum\limits_{\eta\in
V_s:\deg(\eta)=\frac{k}{D(q-1)}}1.$$
\end{prop}

Define
$$\bar{B}_s:=\mathbb{F}_q[x^{L_s(\Delta_{\infty})}]:
=\{\sum\limits_{u\in
L_s(\Delta_{\infty})}a_ux^u:a_u\in\mathbb{F}_q, a_u=0\text{ for
all but finitely many }u\},$$ and
$\bar{B}=\oplus_{s=0}^{q-2}\bar{B}_s$. Then $\bar{B}_s$ is a
$\bar{B}_0$-module, with the multiplication rule
$$x^ux^{u'}=\left\{\begin{array}{ll}
x^{u+u'},&\text{ if } u \text{ and }u'\text{ are cofacial}, \\
0, &\text{ otherwise.}
\end{array}\right. $$
Define $$\text{mod }\pi:B_s\rightarrow\bar{B}_s,\
\sum\limits_{u\in
L_s(\Delta_{\infty})}a_u\pi^{D(q-1)p^{m-1}\deg(u)}x^u\mapsto\sum\limits_{u
\in L_s(\Delta_{\infty})}\bar{a}_ux^u,$$ where $\bar{a}_u=a_u+\pi
O$.
\begin{lemma}The
sequence
$$0\rightarrow B_s\stackrel{\pi}{\rightarrow} B_s
\stackrel{\text{mod }\pi}{\rightarrow}\bar{B}_s\rightarrow0$$ is
exact.\end{lemma} For $j=1,\cdots,n$, we define
$$\bar{\partial}_j: \bar{B}_s\rightarrow\bar{B}_s,\
g\mapsto(x_j\frac{\partial }{\partial x_j}+\overline{_jf})g,$$
where
$$\overline{_jf}=\sum\limits_{i=0}^{m-1}\sum\limits_{k=0}^{m-i-1}
\sum\limits_{\deg(p^{m-i-1}u)=1}u_j a_{iu}^{p^k}x^{p^ku}.$$ By
Lemma 2.14, we have the following lemma. \begin{lemma}For
$j=1,\cdots,n$, the diagram
$$\begin{array}{ccc}
  B_s &\stackrel{\text{mod }\pi}{\rightarrow}& \bar{B}_s \\
 \hat{\partial}_j\downarrow &  & \bar{\partial}_j\downarrow \\
  B_s & \stackrel{\text{mod }\pi}{\rightarrow} & \bar{B}_s \\
\end{array}$$
is commutative.\end{lemma} For $l=0,\cdots,n$, we define
$$\bar{K}_{s,l}=\bigoplus\limits_{1\leq i_1<\cdots<i_l\leq n}
\bar{B}_se_{i_1}\wedge\cdots\wedge e_{i_l},$$ and
$$\bar{K}_{l}=\oplus_{s=0}^{q-2}\bar{K}_{s,l}=\bigoplus\limits_{1\leq i_1<\cdots<i_l\leq n}
\bar{B}e_{i_1}\wedge\cdots\wedge e_{i_l}.$$ For $l=1,\cdots,n$, we
define
$$\bar{\partial}:\bar{K}_{s,l}\rightarrow\bar{K}_{s,l-1},\ ge_{i_1}\wedge\cdots\wedge e_{i_l}
\mapsto \sum\limits_{k=1}^l(-1)^{k-1}\bar{\partial}_{i_k}(g)
e_{i_1}\wedge\cdots\wedge\hat{e}_{i_k}\wedge\cdots\wedge e_{i_l},\
i_1<\cdots<i_l.$$ It is easy to see that the sequence
$$\bar{K}_{s,n}\stackrel{\bar{\partial}}{\rightarrow}\bar{K}_{s,n-1}
\stackrel{\bar{\partial}}{\rightarrow}\cdots\stackrel{\bar{\partial}}{\rightarrow}\bar{K}_{s,0}$$
is a complex.
\begin{prop}The sequence
$$0\rightarrow (K_{s,\bullet},\hat{\partial}) \stackrel{\pi}{\rightarrow}
(K_{s,\bullet},\hat{\partial})\stackrel{\text{mod
}\pi}{\rightarrow}
(\bar{K}_{s,\bullet},\bar{\partial})\rightarrow0$$ is
exact.\end{prop}
\begin{corollary}The $\mathbb{F}_q$-space
$H_l(\bar{K}_{s,\bullet},\bar{\partial})$ is 0 if $l=1,\cdots,n$,
and is of dimension equal to the $O$-rank of
$H_l(K_{s,\bullet},\hat{\partial})$ if $l=0$.\end{corollary} By
that corollary, Proposition 3.1 follows from the following
one.
\begin{prop}The $\mathbb{F}_q$-space
$H_0(\bar{K}_{s,\bullet},\bar{\partial})$ is of dimension
$n!\text{Vol}(\Delta_{\infty})$. Moreover, it has a basis
represented a set $\bar{V}_s$ of homogenous elements such that
$$P_s(t)=\sum\limits_{k=0}t^k\sum\limits_{\eta\in
\bar{V}_s:\deg(\eta)=\frac{k}{D(q-1)}}1.$$
\end{prop}

For $j=1,\cdots,n$, we define
$$\overline{_jf}^0=\sum\limits_{i=0}^{m-1}
\sum\limits_{\deg(p^{m-i-1}u)=1}u_j
a_{iu}^{p^{m-i-1}}x^{p^{m-i-1}u}.$$ For $l=1,\cdots,n$, we
define$$\bar{\partial}^0:\bar{K}_{s,l}\rightarrow\bar{K}_{s,l-1},\
ge_{i_1}\wedge\cdots\wedge e_{i_l}\mapsto
\sum\limits_{k=1}^l(-1)^{k-1}\overline{_{i_k}f}^0g
e_{i_1}\wedge\cdots\wedge\hat{e}_{i_k}\wedge\cdots\wedge e_{i_l},\
i_1<\cdots<i_l.$$ Then
$$\bar{K}_{s,n}\stackrel{\bar{\partial}^0}{\rightarrow}\bar{K}_{s,n-1}
\stackrel{\bar{\partial}^0}{\rightarrow}
\cdots\stackrel{\bar{\partial}^0}{\rightarrow}\bar{K}_{s,0}$$ is a
complex.
\begin{lemma}([LW, Proposition 6.6]) The $\mathbb{F}_q$-space
$H_l(\bar{K}_{\bullet},\bar{\partial}^0)$ is 0 if $l=1,\cdots,n$,
and is of finite dimension if $l=0$.
\end{lemma}
\begin{corollary}The $\mathbb{F}_q$-space
$H_l(\bar{K}_{s,\bullet},\bar{\partial}^0)$ is 0 if
$l=1,\cdots,n$, and is of finite dimension if $l=0$.
\end{corollary}
\begin{corollary}The Poincar\'{e} series of $H_l(\bar{K}_{s,\bullet},\bar{\partial}^0)$
over $\mathbb{F}_q$ is $P_s(t)$. In particular, $P_s(t)$ is a
polynomial with coefficients in the set of natural numbers.
\end{corollary}Modifying the argument of [Ho, 2.10, 2.13], we can prove the
following lemma.\begin{lemma}We have
$P_s(t)\mid_{t=1}=n!\text{Vol}(\Delta_{\infty})$.
\end{lemma}
\begin{prop}Let $V$ be a basis of $H_0(\bar{K}_{s,\bullet},\bar{\partial}^0)$ over $\mathbb{F}_q$ consisting of
homogeneous elements. Then $V$ is also a basis of
$H_l(\bar{K}_{s,\bullet},\bar{\partial})$ over $\mathbb{F}_q$.
\end{prop}
{\it Proof. }First, we show that $\bar{K}_{s,0}$ is generated by
$V$ and $\bar{\partial}(\bar{K}_{s,1})$ over $\mathbb{F}_q$.
Otherwise, among elements of $\bar{K}_{s,0}$ which are not linear
combinations of elements of $V$ and
$\bar{\partial}(\bar{K}_{s,1})$, we choose one of lowest degree.
We may suppose that it is of form $\bar{\partial}^0(\xi)$. Let
$\xi^0$ be the leading term of $\xi$. Then
$\bar{\partial}^0(\xi)-\bar{\partial}(\xi^0)$ is not a linear
combination of elements of $V$ and
$\bar{\partial}(\bar{K}_{s,1})$, and is of lower degree than
$\partial^0(\xi)$. This is a contradiction. Therefore
$\bar{K}_{s,0}$ is generated by $V$ and
$\bar{\partial}(\bar{K}_{s,1})$ over $\mathbb{F}_q$. It remains to
show that $\xi=0$ whenever $\xi$ belongs to
$\bar{\partial}(\bar{K}_{s,1})$ and is a linear combination of
elements of $V$. Otherwise, we may choose one element $\zeta$ of
lowest degree such that $\xi=\bar{\partial}(\zeta)$. Let $\zeta^0$
be the leading term of $\zeta$. Then $\bar{\partial}^0(\zeta^0)$
is a linear combination of elements of $V$ since it is the leading
term of $\bar{\partial}(\zeta)$. So we have
$\bar{\partial}^0(\zeta^0)=0$. By the acyclicity of
$(\bar{K}_{s,\bullet},\bar{\partial}^0)$,
$\zeta^0=\bar{\partial}^0(\eta)$ for some $\eta$. The form
$\zeta-\bar{\partial}(\eta)$ is now of lower degree than $\zeta$,
contradicting to our choice of $\zeta$. This completes the proof
of the proposition.
\section{The Newton polygon}
We prove Theorem 1.2. The second statement follows from the first
and Lemma 3.10. The first statement follows from the following two
lemmas.
\begin{lemma}([LW, Theorem 1.3]) The Newton polygon of $\prod_{\chi}L_{f,\chi}$
with respect to $\text{ord}_q$ lies above the degree-$D(q-1)$
Hodge polygon of $\sum\limits_{s=0}^{q-2}P_s(t)$. Moreover, their
endpoints coincide.
\end{lemma}\begin{lemma}The Newton polygon of
$L_{f,\chi}(t)$ with respect to $\text{ord}_q$ lies above the
degree-$D(q-1)$ Hodge polygon of
$\frac{1}{a}\sum\limits_{i=0}^{a-1}P_{sp^i}(t)$.
\end{lemma}
That lemma follows from Corollary 2.18 and the following two
lemmas.
\begin{lemma}The Newton polygon of
$\det_O(1-\phi_p^at;H_0(K_{s,\bullet},\hat{\partial}))$ with
respect to $\text{ord}_q$ is obtained from the Newton polygon of
$\det_{\mathbb{Z}_p[\pi]}(1-\phi_pt;\oplus_{i=0}^{a-1}H_0(K_{sp^i,\bullet},\hat{\partial}))$
with respect to $\text{ord}_p$ by dividing the ordinates and
abscissas by $a^2$.
\end{lemma}
\begin{lemma}The Newton polygon of
$\det_{\mathbb{Z}_p[\pi]}(1-\phi_pt;\oplus_{i=0}^{a-1}H_0(K_{sp^i,\bullet},\hat{\partial}))$
with respect to $\text{ord}_p$ lies above the degree-$D(q-1)$
Hodge polygon of $a\sum\limits_{i=0}^{a-1}P_{sp^i}(t)$.
\end{lemma}

{\it Proof of Lemma 4.3. }Let $\sigma$ acts on $O[t]$
coefficient-wise. Then
$$\text{det}_O(1-\phi_p^at;H_0(K_{s,\bullet},\hat{\partial}))^{\sigma}=
\text{det}_O(1-\phi_p^at;H_0(K_{sp,\bullet},\hat{\partial})).$$ So
$$\text{det}_{\mathbb{Z}_p[\pi]}(1-\phi_p^at;\oplus_{i=0}^{a-1}H_0(K_{sp^i,\bullet},\hat{\partial}))
=\text{N}_{O/\mathbb{Z}_p[\pi]}\text{det}_O(1-\phi_p^at;\oplus_{i=0}^{a-1}H_0(K_{sp^i,\bullet},\hat{\partial}))
$$$$=\text{det}_O(1-\phi_p^at;\oplus_{i=0}^{a-1}H_0(K_{sp^i,\bullet},\hat{\partial}))^a.$$
Hence
$$\prod_{\zeta^a=1}\text{det}_{\mathbb{Z}_p[\pi]}(1-\zeta\phi_pt;\oplus_{i=0}^{a-1}H_0(K_{sp^i,\bullet},\hat{\partial}))
=\text{det}_{\mathbb{Z}_p[\pi]}(1-\phi_p^at^a;\oplus_{i=0}^{a-1}H_0(K_{sp^i,\bullet},\hat{\partial}))
$$$$=\text{det}_O(1-\phi_p^at^a;\oplus_{i=0}^{a-1}H_0(K_{sp^i,\bullet},\hat{\partial}))^a.$$
Note that
$\text{det}_O(1-\phi_p^at;H_0(K_{s,\bullet},\hat{\partial}))$ and
$\text{det}_O(1-\phi_p^at;H_0(K_{sp,\bullet},\hat{\partial}))$,
being conjugate to each other, share the same Newton polygon.
Lemma 4.3 then follows.

{\it Proof of Lemma 4.4. }By Proposition 3.1, the $O$-module
$H_0(K_{s,\bullet},\hat{\partial})$ has a basis represented a set
$V_s$ of homogenous elements such that
$$P_s(t)=\sum\limits_{k=0}t^k\sum\limits_{\eta\in
V_s:\deg(\eta)=\frac{k}{D(q-1)}}1.$$ For real numbers
$b>\frac{1}{p-1}$ and $c$, write
$$L_s(b,c)=\{\sum_{u\in L_s(\Delta_{\infty})}a_ux^u :
 a_u\in \mathbb{Q}_p[\mu_{q-1},\pi], \text{ ord}_p(a_u)\geq
 b\deg(u)+c\},$$ and denote by $V_s(b,c)$ the subset of
elements of $L_s(b,c)$ which are finite linear combinations of
elements of $V_s$. The space $L_s(b,c)$ is compact with respect to
the topology of coefficient-wise convergence. We claim that, if
$\frac{1}{p-1}<b<\frac{p}{p-1}$, then
$$L_s(b,c)=V_s(b,c)+\sum\limits_{k=1}^n\hat{\partial}_kL_s(b,c+b-\frac{1}{p-1}).$$
In fact, that claim follows from a result of Liu-Wei [LW,
Proposition 8.2], which says that, if
$\frac{1}{p-1}<b<\frac{p}{p-1}$, then
$$L(b,c)=V(b,c)+\sum\limits_{k=1}^n\hat{\partial}_kL(b,c+b-\frac{1}{p-1}),$$
where
 $L(b,c)=\sum\limits_{s=0}^{q-2}L_s(b,c)$, and $V=\sum\limits_{s=0}^{q-2}V_s$.

Let $\lambda\in\mu_{q-1}$ such that $\lambda,\cdots, \lambda^a$ is
a basis of $\mathbb{Z}_p[\mu_{q-1},\pi]$ over $\mathbb{Z}_p[\pi]$.
Then
$$V_{\chi}=\{\lambda^i\xi:i=1,\cdots,a,\xi\in
V_{sp^j},j=0,\cdots,a-1\}$$ is a basis of
$\oplus_{i=0}^{a-1}H_0(K_{sp^i,\bullet},\hat{\partial})$ over
$\mathbb{Z}_p[\pi]$. For each $\xi\in V_{\chi}$, we write
$$\phi_p(\xi)\equiv\sum\limits_{\eta\in
V_{\chi}}c_{\eta,\xi}\eta \ (\text{mod
}\sum\limits_{k=1}^n\hat{\partial}_kB_{\chi}),\
c_{\eta,\xi}\in\mathbb{Z}_p[\pi],$$ where
$B_{\chi}=\sum\limits_{i=0}^{a-1}B_{sp^i}$. Then
$$\text{det}_{\mathbb{Z}_p[\pi]}(1-\phi_pt;\oplus_{i=0}^{a-1}H_0(K_{sp^i,\bullet},\hat{\partial}))
=\det(1-(c_{\eta,\xi})t).$$ By definition, $\phi_p(\xi)$ lies in
the space $L(\frac{p}{p-1})$. So, by our claim, $c_{\eta,\xi}\eta$
lies in every $L(b)$ with $\frac{1}{p-1}<b<\frac{p}{p-1}$. Hence
$\text{ord}_p(c_{\eta,\xi})\geq(b-\frac{1}{p-1})\deg(\eta)$ for
every $\frac{1}{p-1}<b<\frac{p}{p-1}$. It follows that
$\text{ord}_p(c_{\eta,\xi})\geq\deg(\eta)$. Write
$$\det(1-(c_{\eta,\xi})t)=\sum\limits_k\alpha_nt^n.$$ And order
elements of $V_{\chi}$ as $\eta_1,\eta_2,\cdots$ such that
$\deg(\eta_i)\leq\deg(\eta_{i+1})$. Then
$$\text{ord}_p(\alpha_n)\geq\sum\limits_{i\leq n}\deg(\eta_i).$$
It follows that the Newton polygon of the characteristic
polynomial of $(c_{\eta,\xi})$ with respect to $\text{ord}_p$ lies
above the degree-$D(q-1)$ Hodge polygon of
$a\sum\limits_{i=0}^{a-1}P_{sp^i}(t)$. Lemma 4.4 is proved.

\section{The Gauss-Heilbronn sums}
We determine the Newton polygon of $L_{f,\chi}(t)$ when $n=1$ and
$f=\sum\limits_{i=0}^{m-1}V^i([c_ix])$ with $c_0=1$. That is, we
prove Theorem 1.3. By Lemma 4.2, it suffices to prove the
following proposiition.
\begin{prop}The Newton polygon of
$\det_{\mathbb{Z}_p[\pi]}(1-\phi_pt;\oplus_{i=0}^{a-1}H_0(K_{sp^i,\bullet},\hat{\partial}))$
with respect to $\text{ord}_p$ coincides with the polygon with
vertices at $(0,0)$ and the points
$$(a^2n,\frac{a^2(n-1))n}{2p^{m-1}}+\frac{an(d_0+\cdots+d_{a-1})}{p^{m-1}(p-1)}),\ n=1,\cdots,p^{m-1}.$$
\end{prop}

We have $\Delta_{\infty}=[0,p^{m-1}]$, $\deg(p^{m-1})=1$,
$D=p^{m-1}$, and $\overline{_1f}^0=x^{p^{m-1}}$. Write
$$\frac{s}{q-1}=-\sum\limits_{l=0}^{+\infty}s_lp^l,\ 0\leq s_l\leq p-1.$$
Then, for each $l=0,\cdots,a-1$,
$\bar{B}_{sp^{a-l}}=x^{\frac{s_l+s_{l+1}p+\cdots+s_{l+a-1}p^{a-1}}{q-1}}\mathbb{F}_q[x]$,
$$\bar{V}_{sp^{a-l}}=\{x^u:u=\frac{s_l+s_{l+1}p+\cdots+s_{l+a-1}p^{a-1}}{q-1}+i,0\leq i\leq p^{m-1}-1\}$$
represents a basis of
$\bar{B}_{sp^{a-l}}/(\overline{_1f}^0)$ over $\mathbb{F}_q$, and
$$V_{sp^{a-l}}=\{\pi^{(q-1)p^{m-1}u}x^u:u=\frac{s_l+s_{l+1}p+\cdots+s_{l+a-1}p^{a-1}}{q-1}+i,0\leq i\leq p^{m-1}-1\}$$
represents a basis of $H_0(K_{sp^{a-l},\bullet},\hat{\partial})$
over $\mathbb{Z}_p[\mu_{q-1},\pi]$. It follows that
$$U_{sp^{a-l}}=x^{\frac{s_l+s_{l+1}p+\cdots+s_{l+a-1}p^{a-1}}{q-1}}\{(\pi_mx)^i:0\leq i\leq p^{m-1}-1\}$$
represents a basis of
$H_0(K_{sp^{a-l},\bullet},\hat{\partial})\otimes_{\mathbb{Z}_p}\mathbb{Q}_p$
over $\mathbb{Q}_p[\mu_{q-1},\pi]$.

Define $A_l=(a_{ij}^{(l)}\pi_m^{(p-1)i+s_l})_{0\leq i,j\leq
p^{m-1}-1}$ is the matrix of $\phi_p$ from
$H_0(K_{sp^{a-l},\bullet},\hat{\partial})\otimes_{\mathbb{Z}_p}\mathbb{Q}_p$
to
$H_0(K_{sp^{a-l-1},\bullet},\hat{\partial})\otimes_{\mathbb{Z}_p}\mathbb{Q}_p$
with respect to those bases. Write
$$\prod\limits_{i=0}^{m-1}E(\pi_{m-i}c_ix)=\sum\limits_{n=-\infty}^{+\infty}\alpha_n(\pi_mx)^n,\
\alpha_n\in\mathbb{Z}_p[\mu_{q-1},\pi].$$
Then
$$\phi_p(x^{\frac{s_l+s_{l+1}p+\cdots+s_{l+a-1}p^{a-1}}{q-1}}(\pi_m
x)^j)=x^{\frac{s_l+s_{l+1}p+\cdots+s_{l+a-1}p^{a-1}}{q-1}}\sum\limits_{i=0}^{+\infty}
\alpha_{pi-j+s_l}^{\sigma^{-1}}\pi_m^{(p-1)i+s_l}(\pi_m x)^i.$$ It
follows that
$$a_{ij}^{(l)}\equiv \alpha_{pi-j+s_l}^{\sigma^{-1}}
(\text{mod }\pi).$$

Fix a $\lambda\in\mu_{q-1}$ such that $\lambda^1,\cdots,\lambda^a$
is a basis of $\mathbb{Z}_p[\mu_{q-1},\pi]$ over
$\mathbb{Z}_p[\pi]$. Let $A_{ij}^{(l)}$ be the matrix of
$$\mathbb{Z}_p[\mu_{q-1},\pi]\rightarrow\mathbb{Z}_p[\mu_{q-1},\pi],\ x\mapsto a_{ij}^{(l)}x$$
with respect to that basis. Then
$$\{\lambda_i\xi:i=1,\cdots,a,\xi\in
U_{sp^{a-l}}\}$$ represents a basis of
$H_0(K_{sp^{a-l},\bullet},\hat{\partial})\otimes_{\mathbb{Z}_p}\mathbb{Q}_p$
over $\mathbb{Q}_p[\pi]$, and
$A^{(l)}=(A_{ij}^{(l)}\pi_m^{(p-1)i+s_l})_{0\leq i,j\leq
p^{m-1}-1}$ is the matrix of $\phi_p$ from
$H_0(K_{sp^{a-l},\bullet},\hat{\partial})\otimes_{\mathbb{Z}_p}\mathbb{Q}_p$
to
$H_0(K_{sp^{a-l-1},\bullet},\hat{\partial})\otimes_{\mathbb{Z}_p}\mathbb{Q}_p$
with respect to that basis. It follows that
$$\{\lambda_i\xi:i=1,\cdots,a,\xi\in
U_{sp^{a-l},l=0,\cdots,a-1}\}$$ represents a basis of
$\sum\limits_{l=0}^{a-1}H_0(K_{sp^{a-l},\bullet},\hat{\partial})\otimes_{\mathbb{Z}_p}\mathbb{Q}_p$
over $\mathbb{Q}_p[\pi]$, and the matrix of $\phi_p$ with respect
to that basis is
$$\left(%
\begin{array}{ccccc}
  0 & A^{(0)} & 0&\cdots& 0 \\
  0 & 0 & A^{(1)} &\cdots&0 \\
   &  &  & \ddots & \\
  0 &0& \cdots & 0 & A^{(a-2)} \\
  A^{(a-1)} & 0 & 0&\cdots & 0 \\
\end{array}%
\right)$$So Proposition 5.1 follows from  the following
one.\begin{prop}The determinants of the matrices
$(A_{ij}^{(l)})_{0\leq i,j\leq n}$, $0\leq n\leq p^{m-1}-1$, are
$p$-adic units.\end{prop} {\it Proof. }The matrix
$(a_{ij}^{(l)})_{0\leq i,j\leq n}$ defines an endomorphism of
$\mathbb{Z}_p[\mu_{q-1},\pi]^{n+1}$ in a standard way. By the
following proposition, it is an automorphism. The matrix
$(A_{ij}^{(l)})_{0\leq i,j\leq n}$, being a matrix of that
automorphism over $\mathbb{Z}_p[\pi]$, has determinant in
$\mathbb{Z}_p[\pi]^{\times}$. The proposition is proved.
\begin{prop}The determinants of the
matrices $(a_{ij}^{(l)})_{0\leq i,j\leq n}$, $0\leq n\leq
p^{m-1}-1$, are $p$-adic units.\end{prop} {\it Proof. }By Lemma
5.1, it suffices to show that the determinants of the matrices
$(\alpha_{pi-j+s_l})_{0\leq i,j\leq n}$ are $p$-adic units.
 Without loss of generality, we
assume that $c_1=\cdots=c_{m-1}=0$. Write
$$\exp(\pi_mx)
=\sum\limits_{n=-\infty}^{+\infty}b_n(\pi_mx)^n,$$ and
$$E(\pi_mx)\exp(-\pi_mx)
=\sum\limits_{n=-\infty}^{+\infty}\beta_n(\pi_mx)^{pn}.$$ Then
$$\alpha_{pi-j+s_l}=\sum\limits_n\beta_{i-n}b_{pn-j+s_l}.$$
So $$(\alpha_{pi-j+s_l})_{0\leq i,j\leq n}=(\beta_{i-j})_{0\leq
i,j\leq n}\times(b_{pi-j+s_l})_{0\leq i,j\leq n}.$$ Hence it
suffices to show that the determinants of the matrices
$(b_{pi-j+s_l})_{0\leq i,j\leq n}$ are $p$-adic units.

We define the matrices $(e_{ij}^{(l)})_{l\leq i\leq n,0\leq j\leq
n-l}$, $l=0,\cdots,n$, by setting
$$e_{ij}^{(l)}=\frac{p^{l-u}(i-u)!b_{pi-j+s_k-u}}
    {(i-l)!\prod\limits_{v=u}^{l-1}(pv-j+s_k-u)},
$$ where $u$ is the largest integer $\leq\min\{\frac{j-s_k}{p-1}+1,l\}$.
We have
$$(e_{ij}^{(0)})_{0\leq i,j\leq n}=(b_{pi-j+s_k})_{0\leq i,j\leq n}$$ and
$$e_{ij}^{(l+1)}=\left\{%
\begin{array}{ll}
    -e_{ij}^{(l)}+\frac{e_{lj}^{(l)}}{e_{l,j+1}^{(l)}}e_{i,j+1}^{(l)}, & \hbox{if }pl-j+s_k\geq0, \\
    e_{i,j+1}^{(l)}, & \hbox{otherwise.} \\
\end{array}%
\right.    $$ So, denoting by $w$ the integer such that
$p(w-1)\leq n-s_k<pw$, we have
$$\det(e_{ij}^{(l)})_{l\leq i\leq n,0\leq j\leq
n-l}=\det(e_{ij}^{(l+1)}),\ l=0,\cdots,w-1,$$ and
$$\det(e_{ij}^{(l)})=
\frac{p^{l-u}(l-u)!b_{pl-n+l+s_k-u}\det(e_{ij}^{(l+1)})}
    {\prod\limits_{v=u}^{l-1}(pv-n+l+s_k-u)},\ l=w,\cdots,n-1,$$ where $u$ is an integer defined by
    $n-s_k-(p-1)u<l\leq
n-s_k-(p-1)(u-1)$. It follows that
$$\det(b_{pi-j+s_k})_{0\leq i,j\leq n}
=\det(e_{ij}^{(w)})_{w\leq i\leq n,0\leq j\leq n-w},$$
$$\det(e_{ij}^{(w)})=\det(e_{ij}^{(n-s_k-(p-1)(w-1)+1)})\prod\limits_{l=w}^{n-s_k-(p-1)(w-1)}
\frac{p^{l-w}(l-w)!}{(p(l-w)+p-1)!},$$
$$\det(e_{ij}^{(n-s_k-(p-1)u+1)})=\det(e_{ij}^{(n-s_k-(p-1)(u-1)+1)})
\prod\limits_{l=n-s_k-(p-1)u+1}^{n-s_k-(p-1)(u-1)}
\frac{p^{l-u}(l-u)!}{(p(l-u)+p-1)!},\ 0<u<w,$$and
$$\det(e_{ij}^{(n-s_k+1)})=\prod\limits_{l=n-s_k+1}^n
\frac{p^ll!}{(pl+p-1)!}.$$Therefore we have
$$\det(b_{pi-j+s_k})_{0\leq i,j\leq n}=\prod\limits_{l=w}^n
\frac{p^{l-u_l}(l-u_l)!}{(p(l-u_l)+p-1)!},$$where $u_l$ is defined
by $n-s_k-(p-1)u_l<l\leq n-s_k-(p-1)(u_l-1)$. In particular,
$\det(b_{pi-j+s_k})_{0\leq i,j\leq n}$ is a $p$-adic unit. The
proposition is proved.

\end{document}